\documentclass[10pt]{article}
\usepackage{amsmath,amsthm,amsfonts,amssymb,amstext}
\theoremstyle{plain}

\newtheorem{theorem}{Theorem}[section]

\newtheorem{lemma}[theorem]{Lemma}

\begin{document}

%COMMANDS FOR INVARIANT MEASURES CHAPTER

%\newtheorem{theorem}{Theorem}[section]
%\newtheorem{corollary}[theorem]{Corollary}
%\newtheorem{lemma}[theorem]{Lemma}
%\newtheorem{proposition}[theorem]{Proposition}

%\theoremstyle{definition}

%\newtheorem{definition}{Definition}[section]
%\newtheorem{example}{Example}
%\newtheorem{remark}[theorem]{Remark}
%\newtheorem{note}[theorem]{Note}

%END COMMANDS FOR INVARIANT MEASURES CHAPTER

%COMMANDS FOR SPECTRAL GEOMETRY CHAPTER

\newcommand{\Z}{\mathbb Z}
\newcommand{\Q}{\mathbb Q}
\newcommand{\R}{\mathbb R}
\newcommand{\N}{\mathbb N}
\newcommand{\C}{\mathbb C}
\newcommand{\F}{\mathbb F}
\newcommand{\f}{\mathbb f}
\newcommand{\K}{\mathbb K}
\newcommand{\HYP}{\mathbb H}
\newcommand{\germ}{\mathfrak}
\newcommand{\g}{\mathfrak{g}}
\newcommand{\h}{\mathfrak{h}}
\newcommand{\W}{\mathcal{W}}
\newcommand{\T}{\mathcal{T}}
\newcommand{\D}{\mathcal{D}}
\newcommand{\A}{\mathcal{A}}

\newcommand{\Char}{\ensuremath{\operatorname{Char}}}
\newcommand{\Diff}{\ensuremath{\operatorname{Diff}}}
\newcommand{\Ind}{\ensuremath{\operatorname{Ind}}}
\newcommand{\Res}{\ensuremath{\operatorname{Res}}}
\newcommand{\DIM}{\ensuremath{\operatorname{dim}}}
\newcommand{\diag}{\ensuremath{\operatorname{diag}}}
\newcommand{\End}{\ensuremath{\operatorname{End}}}
\newcommand{\Isom}{\operatorname{Isom}}
\newcommand{\SU}{\operatorname{SU}}
\newcommand{\GL}{\operatorname{GL}}
\newcommand{\BC}{\operatorname{BC}}
\newcommand{\SO}{\operatorname{SO}}
\newcommand{\SL}{\operatorname{SL}}
\newcommand{\U}{\operatorname{U}}
\newcommand{\Tr}{\operatorname{Tr}}
\newcommand{\su}{\mathfrak{su}}

\newcommand{\Ad}{\operatorname{Ad}}
\newcommand{\ad}{\operatorname{ad}}
\newcommand{\Spec}{\operatorname{Spec}}
\newcommand{\spec}{\operatorname{spec}}
\newcommand{\Normal}{\operatorname{N}}
\newcommand{\Spin}{\operatorname{Spin}}
\newcommand{\Hol}{\operatorname{Hol}}
\newcommand{\Aut}{\operatorname{Aut}}
\newcommand{\Sim}{\operatorname{Sim}}

%%%%%%%%%%   Title & Absrtact   %%%%%%%%%%%%%%

\title{Measures Invariant under the Geodesic Flow and their Projections}
\author{Craig J. Sutton}
\date{}

\maketitle

\begin{abstract} 
\noindent
Let $S^{n}$ be the $n$-sphere of constant positive curvature.
For $n \geq 2$, we will show that a measure on the unit tangent bundle of $S^{2n}$, which is even and  
invariant under the geodesic flow, is not uniquely determined by its projection to $S^{2n}$.
\end{abstract}



\section{Introduction}

The topological entropy, $h_{T}(g)$, of a Riemannian manifold $(M,g)$ is
a geometric invariant which attempts to capture the complexity of the
geodesic flow. In \cite{KKW} it was shown that for a metric $g$ of
negative sectional curvature the function $h_{T}(g_{\lambda})$ is $C^{1}$,
where $g_{\lambda}$, $-\epsilon < \lambda < \epsilon$, is a $C^{2}$-perturbation of $g$, and an explicit 
formula for the derivative
was obtained. As an application of this formula they established the following
interesting result.

\begin{theorem}\label{Thm:KKW}(\cite[p.21 and 28]{KKW})
Let $M$ be a compact surface and let $\mathcal{R}(M)$ denote the
submanifold of negatively curved $C^{2}$ metrics on $M$ having area
equal to 1. Then $h_{T}: \mathcal{R}(M) \to \R$ has a critical point at
$g_{0}$ if and only if the Lebesgue measure $l_{g_{0}}$ and the Margulis
measure $\mu_{g_{0}}$ with respect to $g_{0}$ have the same projection to $M$; that is,
$l_{g_{0}}$ and $\mu_{g_{0}}$ agree on $\pi^{-1}(\mathcal{B}(M)) = \{ \pi^{-1}(A) : A \in \mathcal{B}(M) \}$, 
where $\mathcal{B}(M)$ is the $\sigma$-algebra of Borel subsets on $M$
and $\pi : S(M) \to M$ is the canonical projection.
\end{theorem}

\noindent
Katok, Knieper and Weiss went on to conjecture that for an arbitrary 
compact manifold $(M,g)$ of negative sectional curvature the Margulis 
and Lebesgue measures with respect to $g$ coincide whenever they have the same
projection to $M$. As they note, establishing this conjecture would then 
demonstrate Theorem~\ref{Thm:KKW} in arbitrary dimensions \cite[p. 21]{KKW}. 

The above results led Flaminio to consider the general problem of determining 
the measures on the unit tangent bundle which are invariant under the geodesic 
flow and are determined by their projection to $M$. By restricting his 
attention to the class of even measures; that is, measures on $S(M)$ which 
are invariant under the flip map $(x, v)\stackrel{\sigma}{\to} (x, -v)$ on 
$S(M)$, Flaminio obtained the following.

\begin{theorem}\label{Thm:Flaminio}(\cite{Flam})
Let $g$ be a metric of positive sectional curvature on $S^{2}$. 
Then an even, $G_{t}$-invariant distribution $T$ is determined by its 
projection to $S^{2}$. That is, $T$ is determined by its values on the set 
$$\pi^{*}(C^{\infty}(S^{2})) \equiv \{ f \circ \pi : f \in C^{\infty}(S^{2}) \},$$ 
where $\pi: S(S^{2}) \rightarrow S^{2}$ is the natural projection.
\end{theorem}

In particular, this result shows that for a closed surface $(M,g)$ of positive
sectional curvature the even, $G_{t}$-invariant probability measures are
determined by their projections to $M$. It is natural to wonder whether this
result generalizes to all closed Riemannian manifolds of positive sectional
curvature. By studying the right regular representation of $SO(n)$ we obtain the
following negative answer. 

\begin{theorem}\label{Thm:InvMeas}
Let $(S^{2j},g)$ be the standard sphere of constant curvature 1 with $j \geq 2$. 
Then even, $G_{t}$-invariant complex measures on $S(S^{2j})$ are not determined by
their projection to $S^{2j}$.
In particular, there are non-zero, even, $G_{t}$-invariant finite real
measures on $S(S^{2j})$ which project to zero on $S^{2j}$.
\end{theorem}

%%%%%%%%%%%%%%%%%%%%%%%%%%%%%%%%%%%%%%%%%%%%%%%%%%%%%%%%%%%%%%%%%%%%%%%%%%%%%%%
\section{Constructing the Measure}

By the Riesz representation theorem there is an isomorphism between
bounded linear functionals and complex meaures. For a linear functional 
$F : L^{2}(S(M)) \to \C$ the notion of projecting to $M$ translates into
restricting $F$ to the set $\pi^{*}(L^{2}(M)) \equiv \{ f \circ \pi : f \in L^{2}(M)
\}$. The notions of evenness and $G_{t}$-invariance are also defined in the obvious
way for $F$. Consequently, we see that we can construct measures as in
Theorem~\ref{Thm:InvMeas} by finding a linear functional with the
corresponding properties.
This will be carried out in the remainder of this paper.

For our discussion we fix the following notation.
\begin{enumerate}
\item[A)] $G=\SO(n)$
\item[B)] $H=\SO(2)\oplus I_{n-2}$
\item[C)] $K=I_{2}\oplus \SO(n-2)$
\item[D)] $L=[1]\oplus \SO(n-1)$
\item[E)] $L^{2}(G,dx) = \{f:G \rightarrow \C$ measureable : 
$\int_{G} {\|f\|}^{2}dx < \infty \}$; where $dx$ is Haar measure.
\item[F)] $G$ will act on $L^{2}(G,dx)$ via the right regular representation 
$\Phi : G \rightarrow \mbox{Aut}(L^{2}(G,dx))$, which is given by $(\Phi(g).f)(x) = f(xg)$.
\item[G)] For any representation $(V, \tau)$ of an arbitrary group $B$ we let 
$V^{S}=\{v \in V : \tau(s).v = v$ \mbox{ for all }$s \in S \}$ for any $S \subset B$.
\end{enumerate}
\noindent
We also note that for the sphere $S^{n-1} = G/L$ of constant positive sectional 
curvature 1 the geodesic flow is given by the right action of $H$ on 
$S(S^{n-1}) = G/K$ and the flip map $\sigma: S(S^{n-1}) \to S(S^{n-1})$ can be 
realized as the right action of $[1]\oplus -I_{n-1}$ on $G/K$ when $n$ is odd. 
%Until further notice we shallrestrict our attention  to $n = 2j + 1, j\geq 2$.

\medskip
In constructing the desired measure we will find the following lemma to be
useful.

\begin{lemma}\label{Lemma:Main}
Let $n = 2j +1 \geq 5$ and let $S$ be the subgroup of
$G$ generated by $H, K \leq G$ and $\sigma \in G$. Then there exists a finite 
dimensional unitary representation $\tau: G \to \GL(W)$ such that 
$(W^{L})^{\perp} \cap W^{S} \neq\{0\}$.
%, where $S = <H,K, \sigma> \leq G$. 
In particular, we may take $(W, \tau)$ to be an irreducible representation of $G$.
\end{lemma}

\noindent
Indeed, let $(W, \tau)$ be an irreducible representation of $G$ as in
Lemma~\ref{Lemma:Main}. Then
$W$ can be thought of as a subrepresentation of $L^{2}(G)$. Using this identification
and taking $\theta \in (W^{L})^{\perp} \cap W^{S}$ we may define 
$\tilde{F} : L^{2}(G) \to \C$ by 
$$f \mapsto \int_{G} f \overline{\theta} dx.$$
It then follows from the Peter-Weyl
theorem (see \cite{Knapp2}) and the way $\theta$ was chosen that
$\tilde{F}$ has the following properties.
\begin{enumerate}
\item\label{Item1} $\tilde{F} |{L^{2}(G)^{L}} \equiv 0.$
\item\label{Item2} $\tilde{F} |{L^{2}(G)^{S}} \not\equiv 0$. In paricular, $\tilde{F}|{L^{2}(G)^{K}}
\not\equiv 0$. 
\end{enumerate}
Now, since $L^{2}(G)^{K}$ can be indentified with $L^{2}(G/K)$ we see
from Property~\ref{Item2} that 
$\tilde{F}$ actually defines a nonzero bounded linear functional 
$F : L^{2}(G/K) \to \C$, which is even and invariant under the geodesic flow.
Furthermore, it follows from Property~\ref{Item1} that the projection of $F$ to $S^{n-1} = G/L$ is zero.
Then as noted earlier the Riesz representation theorem provides us with a non-zero, even,
$G_{t}$-invariant complex measure on $S(S^{n-1})$ that projects to zero. Consequently,
one of the real measures $\mathrm{Re}(\mu)$ or $\mathrm{Im}(\mu)$ also has these properties.
All that remains is to prove is Lemma~\ref{Lemma:Main}.

\begin{proof}[Proof of Lemma \ref{Lemma:Main}]

We let $\germ{g}$, $\germ{l}$, $\germ{h}$, $\germ{k}$ and $\germ{s}$ denote the Lie algebras of
$G$, $L$, $H$, $K$, and $S$ respectively. 
Following an argument due to G. Prasad we will show that 
$\mathrm{Ad} : G \rightarrow \mathrm{Aut}(\mathrm{Sym}^{2}(\germ{g}^{\C}))$, where 
$\germ{g}^{\C} = \germ{g} \oplus i\germ{g}$ is the complexification of $\germ{g}$ and
$\mathrm{Ad}$ is the natural linear extention of the adjoint representation of $G$, is a
representation of $G$ which satisfies Lemma~\ref{Lemma:Main}.  

Upon inspection we can see that 
$\germ{s} = \germ{h} \oplus \germ{k}$ and 
$\germ{g} = \germ{l} \oplus \R^{n-1} = \germ{s} \oplus \R^{2n-2}$. 
From this we can see that
$$\germ{g}^{\C} = \germ{l}^{\C} \oplus (\R^{n-1})^{\C} = \germ{h}^{\C} \oplus \germ{k}^{\C} \oplus
(\R^{2n-2})^{\C}.$$
Hence, 
\begin{eqnarray*}
\mathrm{Sym}^{2}(\germ{g}^{\C}) &=& \mathrm{Sym}^{2}(\germ{l}^{\C}) \oplus\mathrm{Sym}^{2}((\R^{n-1})^{\C}) \oplus 
   (\germ{l}^{\C}\otimes (\R^{n-1})^{\C}) \\
&=& \mathrm{Sym}^{2}(\germ{l}^{\C}) \oplus\mathrm{Sym}^{2}((\R^{n-1})^{\C}) \oplus 
    \mathrm{Hom}( \germ{l}^{\C},(\R^{n-1})^{\C}),
\end{eqnarray*}
and
$$\mbox{Sym}^{2}(\germ{g}^{\C})^{L} = \mbox{Sym}^{2}(\germ{l}^{\C})^{L} 
\oplus \mbox{Sym}^{2}((\R^{n-1})^{\C})^{L} \oplus \mathrm{Hom}( \germ{l}^{\C},(\R^{n-1})^{\C})^{L}.$$
Since $(\mathrm{Ad}_{L},\germ{l}^{\C} )$ and $(\mathrm{Ad}_{L},(\R^{n-1})^{\C})$ are inequivalent irreducible 
representations of $L$ we see that $\mathrm{Hom}( \germ{l}^{\C},(\R^{n-1})^{\C})^{L} = 0$. 
Otherwise we would have a non-zero, $\C$-linear map $T : \germ{l}^{\C} \rightarrow (\R^{n-1})^{\C}$
such that $\mathrm{Ad}_{L}(x) \circ T = T \circ \mathrm{Ad}_{L}(x)$ for all $x \in L$. 
It would then follow from Schur's lemma that $T$  would have to be an isomorphism, which would
contradict the non-equivalence of the representations. Therefore, 
$$\mbox{Sym}^{2}(\germ{g}^{\C})^{L} = \mbox{Sym}^{2}(\germ{l}^{\C})^{L} \oplus \mbox{Sym}^{2}((\R^{n-1})^{\C})^{L}.$$
We now recall the following well-known fact.

\begin{lemma}\label{Lem:Unitary}
Let $G$ be a compact toplogical group and $(\tau, V)$ a finite dimensional irreducible $\C$-representation. Then 
$\DIM_{\C}\mbox{Sym}^{2}(V)^{G} = 1$. That is, up to scalar multiple there is a unique Hermitian inner product 
$\omega \in \mbox{Sym}^{2}(V)$ on $V$ with respect to which $(\tau, V)$ is unitary.
\end{lemma}

Hence, it follows that $\DIM(\mbox{Sym}^{2}(\germ{g}^{\C})^{L}) = 2$. Also, since  
\begin{eqnarray*}
\mbox{Sym}^{2}(\germ{g}^{\C})^{S} & = &  \mbox{Sym}^{2}(\germ{h}^{\C})^{S} \oplus \mbox{Sym}^{2}(\germ{k}^{\C})^{S}
\oplus 
\mathrm{Sym}^{2}((\R^{2n-2})^{\C})^{S} \oplus \mathrm{Hom}( \germ{h}^{\C},\germ{k}^{\C})^{S} \oplus \\
 & &  \mathrm{Hom}( \germ{h}^{\C},(\R^{2n-2})^{\C})^{S} \oplus \mathrm{Hom}( \germ{k}^{\C},(\R^{2n-2})^{\C})^{S}
\end{eqnarray*}
a similar argument shows that $\DIM_{\C}\mathrm{Sym}^{2}(\germ{g}^{\C})^{M} = 3$.
Putting all of this together we see that $\DIM_{\C}(\mbox{Sym}^{2}(\germ{g}^{\C})^{L})^{\perp}) = 
\mbox{Sym}^{2}(\germ{g}^{\C}) -2$ and $\DIM_{\C}\mbox{Sym}^{2}(\germ{g}^{\C})^{S} =3$ , which implies
$$1 \leq \DIM_{\C} (\mbox{Sym}^{2}(\germ{g}^{\C})^{L})^{\perp} \cap \mbox{Sym}^{2}(\germ{g}^{\C})^{S} \leq 3,$$
which proves our lemma.
\end{proof}
 
We point out that our dimension argument fails when $n=3$---as it should by Theorem~\ref{Thm:Flaminio}.
In this case $H=<I_{3}>$ and $L = [1] \oplus SO(2)$.
Hence, $\germ{s} = \germ{k} \oplus \germ{h} = \germ{k} \cong \germ{so}(2) \cong \germ{l}$, which implies that 
$\DIM_{\C}\mbox{Sym}^{2}(\germ{g}^{\C})^{M} = 2 = \DIM_{\C}\mbox{Sym}^{2}(\germ{g}^{\C})^{L}$. 
This prevents the last line of our argument from working. 

%\addcontentsline{toc}{section}{References}
%\bibliography{InvarMeas}

\begin{thebibliography}{KKW91}

\bibitem[Fla92]{Flam}
Livio Flaminio.
\newblock Une remarque sur les distribution invariantes par les flots
  g\'{e}od\'{e}siques des surface.
\newblock {\em C.R. Acad. Sci. Paris, S\'{e}rie I}, 315:735--738, 1992.

\bibitem[KKW91]{KKW}
A.~Katok, G.~Knieper, and H.~Weiss.
\newblock Formulas for the derivative and critical points of topological
  entropy for anosov and geodesic flows.
\newblock {\em Comm. Math. Phys.}, 138:19--31, 1991.

\bibitem[Kna86]{Knapp2}
Anthony~W. Knapp.
\newblock {\em Representation Theory of Semisimple Lie Groups: An Overview
  Based on Examples}.
\newblock Princeton University Press, 1986.

\end{thebibliography}

\bigskip 
\noindent
Craig J. Sutton \\
Department of Mathematics \\
University of Pennsylvania \\
Philadelphia, PA 19104-6395 \\
USA \\
e-mail: cjsutton@math.upenn.edu    
   
\end{document}